\documentclass{amsart}

\usepackage{fancyhdr}
\usepackage{lastpage}
\usepackage{stmaryrd,yhmath}

\newcommand{\bdis}{\begin{displaymath}}
\newcommand{\edis}{\end{displaymath}}
\newcommand{\be}{\begin{equation}}
\newcommand{\ee}{\end{equation}}
\newcommand{\mbb}{\mathbb}
\newcommand{\mcal}{\mathcal}

\newcommand{\vp}{\varphi}

\newcommand{\vth}{\vartheta}

\newcommand{\zf}{\zeta\left(\frac{1}{2}+it\right)}

\theoremstyle{definition}

\theoremstyle{remark}
\newtheorem{remark}[]{Remark}

\newtheorem*{mydef11}{{\bf Theorem 1}}

\newtheorem*{mydef12}{{\bf Theorem 2}}

\newtheorem*{mydef2}{{\bf Definition}}

\newtheorem*{mydef4}{{\bf Corollary}}

\numberwithin{equation}{section}

\begin{document}

\title{Jacob's ladders, $\mcal{Z}_{\zeta,Q^2}$-transformation of real elementary functions and telegraphic
signals generated by the power functions}

\author{Jan Moser}

\address{Department of Mathematical Analysis and Numerical Mathematics, Comenius University, Mlynska Dolina M105, 842 48 Bratislava, SLOVAKIA}

\email{jan.mozer@fmph.uniba.sk}


\begin{abstract}
In this paper we show that the $\mcal{Z}_{\zeta,Q^2}$-transformation of every unbounded signal based on
increasing power function is a telegraphic signal, i.e. the unit rectangular signal.
\end{abstract}
\maketitle

\section{Introduction}

\subsection{}

In this paper we use, instead of \cite{4}, (2.1), the following $(\zeta,Q^2)$-oscillating system
\be \label{1.1}
G^2(x_1,\dots,x_k;y_1,\dots,y_k)=\prod_{r=1}^k \left|\frac{\zeta(1/2+ix_r)}{\zeta(1/2+iy_r)}\right|^2,\
k\leq k_0,\ k_0\in\mbb{N}
\ee
($k_0$ is arbitrary and fixed). Next, we obtain for given admissible elementary real function
\be \label{1.2}
f(t),\ t\in [T,T+U],\ T>0
\ee
the following factorization formula
\bdis
\left|\frac{\zeta(1/2+i\alpha_r)}{\zeta(1/2+i\beta_r)}\right|\sim g[T,U,\alpha_0(T,U,k;f)],\ T\to\infty,
\edis
where $\alpha_0$ obeys
\bdis
0<\alpha_0(T,U)-T<U,
\edis
and (see \cite{4}, (4.7))
\be \label{1.3}
\alpha_0(T,U)=\vp_1^k[d(T,U)],\ d\in (\overset{k}{T},\overset{k}{\wideparen{T+U}}).
\ee
Here, $\vp_1$ is the Jacob's ladder and $\vp_1^k$ is the $k$-th iteration of the $\vp_1$. Finally, we define the
$\mcal{Z}_{\zeta,Q^2}$-transformation of the given function (\ref{1.2}) as follows
\bdis
\begin{pmatrix}
f(t) \\ t\in [T,T+U] \\ U\in (0,U_0)
\end{pmatrix}
\xrightarrow{\mcal{Z}_{\zeta,Q^2}}
\begin{pmatrix}
g[T,U,\alpha_0(T,U,k;f)] \\ U\in (0,U_0) \\ \alpha_0\in (T,T+U]
\end{pmatrix},\ T>T_0>0 ,
\edis
for admissible $U_0>0$. Let us put, for brevity,
\be \label{1.4}
g[T,U,\alpha_0(T,U,\alpha_0(T,U))]=g[U;T],\ U\in (0,U_0),
\ee
for fixed $k$ and $f$.

\subsection{}

The first and the main result concerning the $\mcal{Z}_{\zeta,Q^2}$-transformation is that transformation of
power functions
\bdis
f(t)=t^\Delta,\ t\in [T,T+U],\ 0<U<U_0=o\left(\frac{T}{\ln T}\right),\ \Delta\in\mbb{R}
\edis
($\Delta$ being arbitrary and fixed) is given by
\be \label{1.5}
\begin{pmatrix}
t^\Delta \\ t\in [T,T+U] \\ U\in (0,U_0)
\end{pmatrix}
\xrightarrow{\mcal{Z}_{\zeta,Q^2}}
\begin{pmatrix}
1 \\ U\in (0,U_0) \\ \alpha_0\in (T,T+U)
\end{pmatrix},
\ee
i.e. (see (\ref{1.4}))
\bdis
g(U;T)=1,\ T\to\infty
\edis
for every fixed element of the mentioned class.

\subsection{}

The result (\ref{1.5}) may be interesting from the point of view of the transformations of the
deterministic signals (pulses) in the theory of communication. From this point of view, we can call
the $\mcal{Z}_{\zeta,Q^2}$-transformation as the $\mcal{Z}_{\zeta,Q^2}$-device (comp. (\ref{1.1})).

\begin{remark}
From (\ref{1.5}) it follows, for example, that the unbounded signal
\bdis
\begin{pmatrix}
t^{1000} \\ t\in [L,L+U] \\ U\in (0,1/2)
\end{pmatrix},\quad \forall L>L_0,\ L_0\in\mbb{N}
\edis
on input of the $\mcal{Z}_{\zeta,Q^2}$-device is transformed by this device into telegraphic signal
\bdis
\begin{pmatrix}
1 \\ U\in (0,1/2) \\ \alpha_0\in (L,L+U)
\end{pmatrix},\quad \forall L>L_0,\ L_0\in\mbb{N},
\edis
i.e. into the unit rectangular signal.
\end{remark}

Now, we give, for completeness, the following \emph{opposite} example.
\bdis
\begin{pmatrix}
(t-L)^\Delta \\ t\in [L,L+U] \\ U\in (0,1/2)
\end{pmatrix}
\xrightarrow{\mcal{Z}_{\zeta,Q^2}}
\begin{pmatrix}
\frac{1}{\Delta+1}\left(\frac{U}{\alpha_0-L}\right)^\Delta \\ U\in (0,1/2) \\ \alpha_0\in (L,L+U)
\end{pmatrix},\ L>L_0,\ \Delta>0.
\edis
In this case:
\begin{itemize}
\item[(a)] bounded simply behaved signal is on the input of the $\mcal{Z}_{\zeta,Q^2}$-device,
\item[(b)] on the output of $\mcal{Z}_{\zeta,Q^2}$-device we obtain the signal of completely
indefinite behavior since the distribution of values
\bdis
\alpha_0=\vp_1^k[d(U,T)],\ U\in (0,1/2)
\edis
is unknown.
\end{itemize}

\begin{remark}
Let us notice that to date there is no result about interaction of the real elementary functions with the
Riemann's zeta-function $\zf$.
\end{remark}

\section{Definition of the $\mcal{Z}_{\zeta,Q^2}$-transformation}

\subsection{The first application of the mean-value theorem}

Let
\be \label{2.1}
f(t), \ t\in [T,T+U],\ \forall\- T>T_0[f],\ U>0,
\ee
be continuous. Since
\be \label{2.2}
\int_T^{T+U}f(t){\rm d}t=F(T+U)-F(T)=U\frac{F(T+U)-F(T)}{U}=UH(T,U),
\ee
then
\be \label{2.3}
\frac 1U\int_T^{T+U}f(t){\rm d}t==H(T,U).
\ee

\subsection{The second application of the mean-value theorem}

If (see \cite{3}, (7.1), (7.2))
\be \label{2.4}
0<U<U_0=o\left(\frac{T}{\ln T}\right)
\ee
then we have (comp. \cite{4}, (4.2)--(4.19)) that
\be \label{2.5}
\begin{split}
& \int_T^{T+U}f(t){\rm d}t= \\
& = \int_{\overset{k}{T}}^{\overset{k}{\wideparen{T+U}}}
f[\vp_1^k(t)]\prod_{r=0}^{k-1}\tilde{Z}^2[\vp_1^r(t)]{\rm d}t= \\
& =(\overset{k}{\wideparen{T+U}}-\overset{k}{T})f(\alpha_0)\prod_{r=0}^k \tilde{Z}^2(\alpha_r),
\end{split}
\ee
where (see \cite{4}, (4.7))
\bdis
\begin{split}
& \alpha_r=\vp_1^{k-r}(d),\ r=0,1,\dots,k, \\
& \alpha_r=\alpha_r(T,U,k;f),\ \alpha_0\in (T,T+U).
\end{split}
\edis
Hence, (see (\ref{2.1}), (\ref{2.5}))
\be \label{2.6}
\prod_{r=1}^k \tilde{Z}^2(\alpha_r)=\frac{U}
{\overset{k}{\wideparen{T+U}}-\overset{k}{T}}\frac{H(T,U)}{f(\alpha_0)}.
\ee

\subsection{The third application of the mean-value theorem}

If
\bdis
f(t)=1
\edis
then we have (comp. \cite{4}, (4.16))
\bdis
\prod_{r=1}^k \tilde{Z}^2(\beta_r)=\frac{U}{\overset{k}{\wideparen{T+U}}-\overset{k}{T}},\
\beta_r=\beta_r(T,U,k),\ r=1,\dots,k.
\edis
Finally, we have (comp. \cite{4}, (4.11), (4.13))
\be \label{2.7}
\begin{split}
& \prod_{r=1}^k\left|\frac{\zeta(1/2+i\alpha_r)}{\zeta(1/2+i\beta_r)}\right|^2= \\
& = \left\{ 1+ \mcal{O}\left(\frac{\ln\ln T}{\ln T}\right)\right\}\frac{H(T,U)}{f(\alpha_0)}= \\
& =\left\{ 1+ \mcal{O}\left(\frac{\ln\ln T}{\ln T}\right)\right\}g[T,U,\alpha_0(T,U,k;f)],\ T\to\infty.
\end{split}
\ee

\subsection{}

Motivated by the formulae (\ref{2.3}), (\ref{2.5}) and (\ref{2.7}), we give the following

\begin{mydef2}
We define the $\mcal{Z}_{\zeta,Q^2}$-transformation acting on the subset $\{ f(t)\}$ of the class of real elementary functions as follows
\be \label{2.10}
\begin{pmatrix}
f(t) \\ t\in [T,T+U] \\ U\in (0,U_0)
\end{pmatrix}
\xrightarrow{\mcal{Z}_{\zeta,Q^2}}
\begin{pmatrix}
g[T,U,\alpha_0(T,U,k;f)] \\ U\in (0,U_0) \\ \alpha_0\in (T,T+U)
\end{pmatrix}
\ee
for every fixed
\bdis
T>\bar{T}_0=\bar{T}_0[f,\vp_1]\geq \bar{T}_0[f].
\edis
\end{mydef2}

\begin{remark}
We put, for brevity, that
\be \label{2.11}
g[T,U,\alpha(T,U,k;f)]=g(U;T)
\ee
for any admissible $k$ and $f$.
\end{remark}

\section{The $\mcal{Z}_{\zeta,Q^2}$-transformation of the power functions}

Let
\be \label{3.1}
\begin{split}
& f(t)=f(t;\Delta)=t^\Delta,\ t\in [T,T+U],\ T\geq \bar{T}_0, \\
& \Delta\in\mbb{R},\ 0<U<U_0=o\left(\frac{T}{\ln T}\right).
\end{split}
\ee

\subsection{}

In the case
\bdis
\Delta\not= -1,0
\edis
we have
\bdis
\begin{split}
& \int_T^{T+U} t^\Delta{\rm d}t=\frac{1}{\Delta+1}\left\{ (T+U)^\Delta-T^\Delta\right\}= \\
& =\frac{1}{\Delta+1}T^{\Delta+1}\left\{\left( 1+\frac{U}{T}\right)^{\Delta+1}-1\right\}= \\
& =\frac{1}{\Delta+1}T^{\Delta+1}\left\{ (\Delta+1)\frac UT+\mcal{O}\left((\Delta+1)^2\frac{U^2}{T^2}\right)\right\}= \\
& = UT^\Delta\left\{ 1+\mcal{O}\left(\frac UT\right)\right\}= \\
& =UT^\Delta\left\{ 1+\mcal{O}\left(\frac{1}{\ln T}\right)\right\},
\end{split}
\edis
i.e.
\be \label{3.2}
\frac 1U\int_T^{T+U} t^\Delta{\rm d}t=T^{\Delta}\left\{ 1+\mcal{O}\left(\frac{1}{\ln T}\right)\right\}.
\ee
Now, from (\ref{3.2}) by (\ref{2.7}) the formula
\be \label{3.3}
\begin{split}
& \prod_{r=1}^k\left|\frac{\zeta(1/2+i\alpha_r^1)}{\zeta(1/2+i\beta_r^1)}\right|=
\left\{ 1+\mcal{O}\left(\frac{\ln\ln T}{\ln T}\right)\right\}\left(\frac{T}{\alpha_0^1}\right)^\Delta, \\
& \alpha_r^1(T,U,k;\Delta),\ r=0,1,\dots,k, \\
& \beta_r^1=\beta_r^1(T,U,k),\ r=1,\dots,k, \\
& \alpha_0^1\in (T,T+U),\ \Delta\not=-1,0,
\end{split}
\ee
follows. Since (see (\ref{3.1}))
\be \label{3.4}
\begin{split}
 & \frac{T}{\alpha_0^1}=\frac{T}{T+\alpha_0^1-T}=\frac{1}{1+\frac{\alpha_0^1-T}{\alpha_0^1}}=1+\mcal{O}\left(\frac{U}{T}\right)=\\
 & =1+\mcal{O}\left(\frac{1}{\ln T}\right),
\end{split}
\ee
then we have (see (\ref{3.3}), (\ref{3.4})) the following formula
\be \label{3.5}
\prod_{r=1}^k\left|\frac{\zeta(1/2+i\alpha_r^1)}{\zeta(1/2+i\beta_r^1)}\right|^2=1+\mcal{O}\left(\frac{\ln\ln T}{\ln T}\right),\
\Delta\not=-1,0.
\ee

\subsection{}

In the case
\bdis
\Delta=-1
\edis
we obtain
\bdis
\frac 1U\int_T^{T+U}\frac{{\rm d}t}{t}=\frac 1T\left\{ 1+\mcal{O}\left(\frac{1}{\ln T}\right)\right\},
\edis
and consequently that
\be \label{3.6}
\prod_{r=1}^k\left|\frac{\zeta(1/2+i\alpha_r^1)}{\zeta(1/2+i\beta_r^1)}\right|^2=1+\mcal{O}\left(\frac{\ln\ln T}{\ln T}\right).
\ee

\subsection{}

Since
\bdis
\Delta=0 \ \Rightarrow \ f(t)=1
\edis
we have directly (comp. (\ref{2.7}) the formula
\be \label{3.7}
\prod_{r=1}^k\left|\frac{\zeta(1/2+i\alpha_r^1)}{\zeta(1/2+i\beta_r^1)}\right|^2=1+\mcal{O}\left(\frac{\ln\ln T}{\ln T}\right).
\ee

\subsection{}

Consequently, we have obtained (see (\ref{2.10}), (\ref{3.5})--(\ref{3.7})) the following

\begin{mydef11}
Under the assumptions (\ref{3.1}) the following holds true
\bdis
\begin{pmatrix}
t^\Delta \\ t\in [T,T+U] \\ U\in (0,U_0)
\end{pmatrix}
\xrightarrow{\mcal{Z}_{\zeta,Q^2}}
\begin{pmatrix}
1 \\ U\in (0,U_0) \\ \alpha_0\in (T,T+U)
\end{pmatrix}, \ T\geq \bar{T}_0[\Delta,\vp_1].
\edis
\end{mydef11}

\section{The $\mcal{Z}_{\zeta,Q^2}$-transformation of unbounded and negligible signals into
telegraphic signals}

\subsection{}

Let
\bdis
\begin{pmatrix}
t^\Delta \\ t\in [L,L+U] \\ U\in (0,a_L) \\ a_L\in (0,1/2]
\end{pmatrix},\
\Delta\not=0,\ L>\bar{L}_0=\bar{L}_0[\Delta,\vp_1],\ L\in\mbb{N}
\edis
be the signal on the input of the $\mcal{Z}_{\zeta,Q^2}$-device (see (\ref{1.3})). Then we have the following

\begin{mydef4}
\be \label{4.1}
\begin{pmatrix}
t^\Delta \\ t\in [L,L+U] \\ U\in (0,a_L) \\ a_L\in (0,1/2]
\end{pmatrix}\xrightarrow{\mcal{Z}_{\zeta,Q^2}}
\begin{pmatrix}
1 \\ U\in (0,a_L) \\ \alpha_0^1\in (L,L+U)
\end{pmatrix},\quad \forall\- L>\bar{L}_0,\ \Delta\not=0,
\ee
i.e. in this case we have (comp. (\ref{2.1})) that
\bdis
g[L,U,\alpha_0^1(L,U,k;\Delta)]=g_L(U;\Delta)=1.
\edis
\end{mydef4}

\begin{remark}
Since
\bdis
\lim_{L\to+\infty}L^\Delta=
\left\{\begin{array}{rcl} +\infty & , & \Delta>0, \\ 0 & , & \Delta<0 ,  \end{array} \right.
\edis
we call the signal (\ref{4.1}) unbounded ($\Delta>0$) and we call it negligible in the case $\Delta<0$.
The $\mcal{Z}_{\zeta,Q^2}$-device transforms those signals into telegraphic one.
\end{remark}

\subsection{}

We shall call the telegraphic signal (comp. (\ref{4.1}))
\be \label{4.2}
S_L(U;a)=
\begin{pmatrix}
1 \\ U\in (0,a) \\ a\in (0,1/2] \\ \alpha_0\in (L,L+a)
\end{pmatrix}
\ee
as the periodic one since
\bdis
S_L(U;a)=S_{L+1}(U;a),\ \forall\- L>\bar{L}_0.
\edis
Next, if the sequence
\bdis
\{ a_L\}_{L>\bar{L}_0}
\edis
is not a stationary one then we shall call the corresponding telegraphic signal as the aperiodic one.

\begin{remark}
We see (comp. (\ref{4.1}) and (\ref{4.2})) that the $\mcal{Z}_{\zeta,Q^2}$-device generates the continuum set of
the periodic signals.
\end{remark}

\begin{remark}
Of course, we may use in (\ref{4.1}) instead of the sequence
\bdis
L,L+1,L+2,\dots ,\ L>\bar{L}_0,\ L\in\mbb{N}
\edis
an arbitrary sequence
\bdis
\{ L_n\}_{n=1}^\infty,\ L_n\in\mbb{R},\ L_1>\bar{L}_0 ,
\edis
and we may put
\bdis
a_n\in \left(\left. 0,\frac{L_{n+1}-L_n}{2}\right.\right];\ U\in (0,a_n).
\edis
\end{remark}

\section{An opposite case: complete uncertainty of a signal on the output of
$\mcal{Z}_{\zeta,Q^2}$-device}

\subsection{}

If
\be \label{5.1}
\begin{split}
& f(t,\Delta)=(t-L)^\Delta,\ t\in [L,L+U],\ \Delta>0,\\
& U\in (0,a_L),\ a_L\in (0,1/2],
\end{split}
\ee
then
\bdis
\int_L^{L+U}(t-L)^\Delta{\rm d}t=\frac{1}{\Delta+1}U^{\Delta+1},
\edis
i.e.
\bdis
\frac 1U \int_L^{L+U}(t-L)^\Delta{\rm d}t=\frac{1}{\Delta+1}U^\Delta,
\edis
and consequently, we have (see (\ref{2.10}), (\ref{4.1})) the following

\begin{mydef12}
\be \label{5.2}
\begin{pmatrix}
(t-L)^\Delta \\ t\in [L,L+U] \\ U\in (0,a_L) \\ a_L\in (0,1/2]
\end{pmatrix}\xrightarrow{\mcal{Z}_{\zeta,Q^2}}
\begin{pmatrix}
\frac{1}{\Delta+1}\left(\frac{U}{\alpha_0^2-L}\right)^\Delta \\ U\in (0,a_L) \\ \alpha_0^2\in (L,L+U)
\end{pmatrix}
\ee
for all
\bdis
L>\bar{L}_0,\ \Delta>0.
\edis
\end{mydef12}

\begin{remark}
In this case the signal
\bdis
\frac{1}{\Delta+1}\left(\frac{U}{\alpha_0^2-L}\right)^\Delta
\edis
on the output of the $\mcal{Z}_{\zeta,Q^2}$-device is completely uncertain since we know about
the values of the function
\bdis
\alpha_0^2=\alpha_0^2(U),\ U\in (0,a_L)
\edis
(for fixed values of $L,k,\Delta$) only that
\bdis
\alpha_0^2(U)-L\in (0,U).
\edis
\end{remark}

\subsection{}

Here we give some remarks about the uncertainty mentioned above. Since for every
\bdis
U\in (0,a_L) \ \Rightarrow \ 0<\alpha_0^2(U)-L<U,\ \Delta>0
\edis
(for fixed $L,k,\Delta$) then either
\be \label{5.3}
0<\alpha_0^2(U)-L\leq \frac U2,
\ee
or
\be \label{5.4}
\frac U2 <\alpha_0^2(U)-L.
\ee
Now, in the case (\ref{5.3}) we have
\be \label{5.5}
\frac{1}{\Delta+1}\left(\frac{U}{\alpha_0^2-L}\right)^\Delta\geq \frac{2^\Delta}{\Delta+1},
\ee
and in the case (\ref{5.4}) we have
\be \label{5.6}
\frac{1}{\Delta+1}\leq  \frac{1}{\Delta+1}\left(\frac{U}{\alpha_0^2-L}\right)^\Delta<\frac{2^\Delta}{\Delta+1}.
\ee

\begin{remark}
Inequalities (\ref{5.5}) and (\ref{5.6}) give some characterization of the distribution of the values
\bdis
\alpha_0^2(U),\ U\in (0,a_L).
\edis
\end{remark}

\section{Properties of the $\mcal{Z}_{\zeta,Q^2}$-transformation}

\subsection*{(A)}

The sequence
\be \label{6.1}
\{ \overset{k}{T}\}_{k=1}^{k_0},\ k_0\in\mbb{N},\ T>T_0>0
\ee
is defined by the formula (comp. \cite{3}, (5.2))
\bdis
\vp_1(\overset{k}{T})=\overset{k-1}{T},\ k=1,\dots,k_0,\ \overset{0}{T}=T.
\edis
Since
\bdis
\overset{k}{T}=\vp_1^{-1}(\overset{k-1}{T})=\vp_1^{-1}(\vp_1^{-1}(\dots \vp_1^{-1}(T)))=\vp_1^{-k}(T),
\edis
then we call the sequence (\ref{6.1}) as the reversely iterated one. Now, we have (similarly to (\ref{6.1})) the following.
To every segment
\bdis
[T,T+U],\ T>T_0
\edis
there is the sequence
\bdis
\{[\overset{r}{T},\overset{r}{\wideparen{T+U}}]\}_{r=0}^{k},\ k\leq k_0
\edis
of the reversely iterated segments, where
\bdis
[\overset{0}{T},\overset{0}{\wideparen{T+U}}]=[T,T+U].
\edis
Since (comp. \cite{4}, (4.7))
\bdis
d\in (\overset{k}{T},\overset{k}{\wideparen{T+U}}) \
\Rightarrow\ \vp_1^r(d)\in (\overset{k-r}{T},\overset{k-r}{\wideparen{T+U}}),\ r=0,1,\dots,k
\edis
then we put
\bdis
\alpha_{k-r}=\vp_1^r(d),\ r=0,1,\dots,k.
\edis
Of course,
\bdis
\alpha_0=\vp_1^k(d).
\edis
For
\bdis
\beta_r, \ r=1,\dots,k
\edis
similar properties hold true (comp. \cite{4}, (4.17)).

\subsection*{(B)}

Next, the sequences
\be \label{6.2}
\{\alpha_r\}_{r=0}^k,\ \{\beta_r\}_{r=0}^k
\ee
have the following properties
\be \label{6.3}
\begin{split}
& T<\alpha_0<\alpha_1<\dots<\alpha_k \\
& T<\beta_1<\dots<\beta_k \\
& \alpha_0\in (T,T+U),\\
& \alpha_r,\beta_r\in (\overset{r}{T},\overset{r}{\wideparen{T+U}}),\ r=1,\dots,k,
\end{split}
\ee
and (see \cite{3}, (5.12))
\be \begin{split}
& \alpha_{r+1}-\alpha_r\sim (1-c)\pi(T),\ r=0,1,\dots,k-1, \\
& \beta_{r+1}-\beta_r\sim (1-c)\pi(T),\ r=1,\dots,k-1,
\end{split}
\ee
where
\be \label{6.5}
\pi(T)\sim \frac{T}{\ln T},\ T\to\infty
\ee
is the prime-counting function and $c$ is the Euler's constant.

\begin{remark}
Jacob's ladder can be viewed, by the formula (see \cite{1}, (6.2))
\bdis
T-\vp_1(T)\sim (1-c)\pi(T),
\edis
as an asymptotic complementary function to the function
\bdis
(1-c)\pi(T)
\edis
in the following sense
\bdis
\vp_1(T)+(1-c)\pi(T)\sim T,\ T\to\infty.
\edis
\end{remark}

\begin{remark}
The asymptotic behavior of the sequences (\ref{6.2}) is as follows: if $T\to\infty$ then the points of every sequence in
(\ref{6.2}) recede unboundedly each from other and all together recede to infinity. Hence, at $T\to\infty$ each sequence
in (\ref{6.2}) behaves as one-dimensional Friedmann-Hubble universe.
\end{remark}

\subsection*{(C)}

Let us remind that the Jacob's ladder
\bdis
\vp_1(t)=\frac 12\vp(t)
\edis
has been introduced in our work \cite{1} (see also \cite{2}), where the function $\vp(t)$ is an arbitrary solution of the
nonlinear integral equation
\bdis
\int_0^{\mu[x(T)}Z^2(t)e^{-\frac{2}{x(T)}t}{\rm d}t=\int_0^T Z^2(t){\rm d}t,
\edis
where each admissible function
\bdis
\mu (y)
\edis
generates the solution
\bdis
y=\vp(T;\mu)=\vp(T),\ \mu(y)\geq ty\ln y.
\edis
The function $\vp_1(T)$ is called the Jacob's ladder corresponding to the Jacob's dream in Chumash, Bereishis, 28:12.

\begin{remark}
By making use of those Jacob's ladders we have shown (see \cite{1}) that the classical Hardy-Littlewood integral (1918)
\bdis
\int_0^T\left|\zf\right|^2{\rm d}t
\edis
has -- in addition to the previously known Hardy-Littlewood expression (and other similar to that one) possessing an unbounded error at
$T\to\infty$ -- the following set of almost exact representations
\bdis
\begin{split}
& \int_0^T\left|\zf\right|^2{\rm d}t=\vp_1(T)\ln T+(c-\ln 2\pi)\vp_1(T)+c_0+\mcal{O}\left(\frac{\ln T}{T}\right),\ T\to \infty ,
\end{split}
\edis
where $c$ is the Euler's constant and $c_0$ is the constant from the Titchmarsh-Kober-Atkinson formula.
\end{remark}

\subsection*{(D)}

We call the system $(\zeta,Q^2)$ as the oscillating system. This is based on the spectral form of the Riemann-Siegel formula (see \cite{4}, (3.8))
\bdis
\begin{split}
& Z(t)=\sum_{n\leq \tau(x_r)}\frac{2}{\sqrt{n}}\cos\{ t\omega_n(x_r)+\psi(x_r)\}+R(x_r), \\
& \tau(x_r)=\sqrt{\frac{x_r}{2\pi}},\ R(x_r)=\mcal{O}(x_r^{-1/4}), \\
& t\in [x_r,x_r+V],\ V\in (0,x_r^{1/4}),
\end{split}
\edis
where the functions
\bdis
\frac{2}{\sqrt{n}}\cos\{ t\omega_n(x_r)+\psi(x_r)\}
\edis
are the Riemann's oscillators with:
\begin{itemize}
\item[(a)] the amplitude
\bdis
\frac{2}{\sqrt{n}},
\edis
\item[(b)] the incoherent local phase constant
\bdis
\psi(x_r)=-\frac{x_r}{2}-\frac{\pi}{8},
\edis
\item[(c)] the nonsynchronized local time
\bdis
t(x_r)\in [x_r,x_r+V],
\edis
\item[(d)] the local spectrum of cyclic frequencies
\bdis
\{ \omega_n(x_r)\}_{n\leq \tau(x_r)},\ \omega_n(x_r)=\frac{\tau(x_r)}{n},
\edis
and similar formulae take place also for $x_r\longrightarrow y_r$.
\end{itemize}

Of course,
\bdis
Z(t)=e^{i\vth(t)}\zf \ \Rightarrow \ |Z(t)|=\left|\zf\right|.
\edis

\thanks{I would like to thank Michal Demetrian for his help with electronic version of this paper.}

\end{document}